\documentclass[oneside,reqno,12pt]{amsart}
\usepackage{amsfonts,amssymb,eucal,amsmath}

\newtheorem{theorem}{Theorem}[section]

\newtheorem{lemma}[theorem]{Lemma}
\newtheorem{proposition}[theorem]{Proposition}
\theoremstyle{definition}
\newtheorem{definition}[theorem]{Definition}
\theoremstyle{remark}
\newtheorem{remark}[theorem]{Remark}
\theoremstyle{remark}

\numberwithin{equation}{section}
\newcommand{\ve}{\varepsilon}
\newcommand{\set}[1]{\left\{#1\right\}}

\renewcommand{\le}{\leqslant}
\renewcommand{\ge}{\geqslant}
\newcommand{\vv}[1]{\mathbf{#1}}
\newcommand{\tvv}[1]{\tilde{\mathbf{#1}}}

\newcommand{\diam}{\operatorname{diam}}
\newcommand{\dist}{\operatorname{dist}}
\newcommand{\NN}{\mathbb{N}}
\newcommand{\ZZ}{\mathbb{Z}}

\newcommand{\RR}{\mathbb{R}}

\newcommand{\B}{\mathcal{B}}

\newcommand{\G}{\mathcal{G}}

\renewcommand{\L}{\mathcal{L}}
\newcommand{\M}{{\mathcal{M}}}
\newcommand{\R}{\mathcal{R}}
\newcommand{\U}{U}

\newcommand{\va}{\mathbf{a}}
\newcommand{\vb}{\mathbf{b}}
\newcommand{\vf}{\mathbf{f}}

\newcommand{\vx}{\mathbf{x}}
\newcommand{\vy}{\mathbf{y}}
\newcommand{\vz}{\mathbf{z}}

\newcommand{\Veronese}{\mathcal{V}}
\newcommand{\ZZnn}{\ZZ_{\ge0}}
\newcounter{Problem}
\newcommand{\Problem}{\refstepcounter{Problem}{\bf
Problem\,\arabic{Problem}.}\,}
\begin{document}


\title[The Khintchine--Groshev theorem on manifolds]
{Metric Diophantine approximation: The Khintchine--Groshev theorem
for non-degenerate manifolds}

\author{V.V.~Beresnevich}
\address{Victor Beresnevich: Department of Number Theory,
Institute of Mathematics, The National Belarus Academy of
Sciences, 220072, Surganova 11, Minsk, Belarus}
\email{beresnevich@im.bas-net.by}

\author{V.I.~Bernik}
\address{Vasily Bernik: Department of Number Theory,
Institute of Mathematics, The National Belarus Academy of
Sciences, 220072, Surganova 11, Minsk, Belarus}
\email{bernik@im.bas-net.by}

\author{D.Y.~Kleinbock}
\address{Dmitry Kleinbock:
Department of Mathematics, Brandeis University, Waltham, MA 02454,
USA} \email{kleinboc@brandeis.edu}

\author{G.A.~Margulis}
\address{Gregory Margulis:
Department of Mathematics, Yale University, New Haven, Connecticut
06520, USA} \email{margulis@math.yale.edu}

\date{January 30, 2002}



\thanks{{\it Thanks}\/: The work of the first and second named authors was
supported in part by grant F00-249 of the Belarussian Fond of
Fundamental Research. The work of the third named author was
supported in part by NSF Grant DMS-0072565, and that of the fourth
named author by NSF Grant DMS-9800607.}


\begin{abstract}
The main objective of this paper is to prove a Khintchine type
theorem for divergence for linear Diophantine approximation on
non-degenerate manifolds, which completes earlier results for
convergence.
\end{abstract}

\maketitle


\section{Background and the main result}
\label{background}

\subsection{Notation}
\label{notation}

The Vinogradov symbol $\ll$ ($\gg$) means ``$\le$ ($\ge$) up to a
positive constant multiplier''; $a\asymp b$ is equivalent to $a\ll
b\ll a$. The usual inner product in $\RR^n$ of $\va$ and $\vb$
will be denoted by $\va\cdot\vb$; $\|\va\|=\sqrt{\va\cdot\va}$ is
the Euclidean norm of $\va$. Also, $\|\va\|_\infty=\max_{1\le i\le
n}|a_i|$ and $\|\va\|_1=\sum_{i=1}^n|a_i|$, where $a_i$ are the
coordinates of $\va$ in the standard basis of $\RR^n$. The
Lebesgue measure of $A\subset\RR^d$ is denoted by $|A|_d$. We
write $|A|$ instead of $|A|_d$ if there is no risk of confusion.
Given a subset $A$ of $\RR^n$, we define
$\diam(A)=\sup_{\va,\vb\in A}\|\va-\vb\|$. Given two subsets $A$
and $B$ of $\RR^n$, we define $\dist(A,B)=\inf_{\va\in A,\vb\in
B}\|\va-\vb\|$; also $\dist(\va,A)=\dist(\{\va\},A)$. Given an
$\vx\in\RR^n$, there is a unique point $\va\in\ZZ^n$ such that
$\vx-\va\in(-1/2,1/2]^n$. This difference will be denoted by
$\langle \vx\rangle$. Given a set $A\subset\RR^d$ and a number
$r>0$, let $\B(A,r)=\{\vx\in\RR^d:\dist(\vx,A)<r\}$. In
particular, $\B(\va,r)= \B(\{\va\},r)$ is the open ball in $\RR^d$
of radius $r$ centered at $\va$. Given a ball $\B=\B(\vx,r)$ and a
positive number $\lambda$, $\lambda\B$ will denote the ball
$\B(\vx,\lambda r)$. Given a map $\vf:\U\longrightarrow\RR^n$, where $\U$
is an open subset of $\RR^d$, we will denote by $\partial_i
\vf:U\longrightarrow \RR^n$, $i = \overline{1,d}$, its partial derivative
with respect to $x_i$. Also we define a map $\nabla \vf:\U\longrightarrow
M_{n\times d}(\RR)$, where $M_{n\times d}(\RR)$ is the space of
$n\times d$ matrices over $\RR$, by setting $\nabla
\vf(\vx)=(\partial_j f_i(\vx))_{1\le i\le n,1\le j\le d}$. We will
also need higher order differentiation: for a {\it multiindex\/}
$\beta = (i_1,\dots,i_d)$, $i_j\in \ZZnn$, where
$\ZZnn=\{x\in\ZZ:x\ge0\}$, we let $\partial_\beta =
\partial_1^{i_1}\circ\dots\circ \partial_d^{i_d}$. Throughout
the paper, $\psi:\RR_+\longrightarrow\RR_+$ is a non-increasing function
unless a different condition is assumed.

\subsection{Metric Diophantine approximation in $\RR^n$.}
\label{MDA in R^n}

Metric Diophantine approximation began with the works of E.\,Borel
and A.J.\,Khintchine, who considered approximation to real numbers
by rational numbers. In 1924 for $n=1$ Khintchine
\cite{Khintchine2} and in 1938 for $n>1$ A.V.\,Groshev
\cite{Groshev1} established a criterion for the solubility of the
inequality
\begin{equation}\label{bbkm:001}
|\,\langle\va\cdot\vy\rangle\,|<\psi(\|\vv a\|_\infty^n)
\end{equation}
in $\va\in\ZZ^n$ for generic $\vy\in \RR^n$. At this point we need
the following

\begin{definition}
The point $\vy\in\RR^n$ is called {\it $\psi$-approximable}\/ if
(\ref{bbkm:001}) has infinitely many solutions $\va\in\ZZ^n$. The
point $\vy\in\RR^n$ is called {\em very well approximable\/} (VWA)
if it is $\psi_\ve$-approximable for some positive $\ve$, where
$\psi_\ve(h)= h^{-(1+\ve)}$.
\end{definition}

In view of this definition, the Khintchine--Groshev theorem
\cite{Khintchine2,Groshev1} asserts that if the sum
\begin{equation}\label{bbkm:002}
\sum_{h=1}^\infty \psi(h)
\end{equation}
diverges (converges), then almost all (almost no) points $\vy\in
\RR^n$ are $\psi$-approximable.

\begin{remark}\label{rem1}
Originally the inequality $|\,\langle\va\cdot\vv
x\rangle\,|<\psi(\|\vv a\|_\infty)$ was considered instead of
(\ref{bbkm:001}). In this setting $\sum_{q=1}^\infty
q^{n-1}\psi(q)$ should be used instead of (\ref{bbkm:002}).
Khintchine assumed that $h\psi(h)$ was non-increasing, and
Groshev's requirement was the monotonicity of $h^{n-1}\psi(h)$.
Later W.M.~Schmidt succeeded to avoid the monotonicity restriction
when $n>1$ (see Section~\ref{conc}).
\end{remark}

\begin{remark}
The Khintchine--Groshev theorem implies that almost all
$\vy\in\RR^n$ are not VWA. The convergence case of the theorem can
be easily derived from the Borel--Cantelli lemma. The main
difficulty is contained in the divergence case.
\end{remark}

\subsection{The concept of Diophantine approximation on
manifolds.}

This concept emerges if one restricts the point $\vy$ to lie on a
submanifold $\M$ of $\RR^n$. Since the manifold $\M$ of dimension
$<n$ itself has zero measure, the Khintchine--Groshev theorem does
not even guarantee the existence of a single $\psi$-approximable
point on $\M$. To make the theory rich in content one tries to
establish if a given property holds for almost all points of this
manifold with respect to the Lebesgue measure induced on the
manifold. We will use the following terminology (more details can
be found in \cite{BernikDodson1} and Section~\ref{conc}).

\begin{definition}\label{defn2}
Let $\M$ be a submanifold of $\RR^n$. One says that $\M$ is {\em
extremal}\/ if almost all points of $\M$ are not VWA. One says
that $\M$ is of {\em Groshev type for divergence $($for
convergence$)$}\/ if almost all (almost no) points of $\M$ are
$\psi$-approximable whenever the sum (\ref{bbkm:002}) diverges
(converges).
\end{definition}

\subsection{Diophantine approximation on the Veronese curves}
\label{Veronese}

In 1932 K.~Mahler \cite{Mahler1} made a conjecture which in the
terminology of this paper claimed that for any $n\in\NN$ the {\em
Veronese curve\/}
\begin{equation}\label{bbkm:003}
\Veronese_n=\{(x,x^2,\dots,x^n): x\in \RR\}
\end{equation}
was extremal. It arose in transcendental number theory in
connection with a classification of real numbers suggested by
Mahler himself. A great deal of work had been undertaken to prove
Mahler's conjecture by J.\,Kubilius, B.\,Volkmann, W.\,LeVeque,
F.\,Kash and W.M.\,Schmidt. In particular, the problem was solved
for $n=2$ by Kubilius \cite{Kubilius1} and for $n=3$ by Volkmann
\cite{Volkmann2}. The complete solution was given by
V.G.~Sprind\^zuk \cite{Sprindzuk1} in 1964.

In 1966 A.~Baker \cite{Baker3} improved Sprind\^zuk's result by
replacing the ``powering'' error function with a general monotonic
function $\psi$ by showing that if
\begin{equation}\label{bbkm:004}
\sum_{k = 1}^\infty \frac{\psi(k)^{1/n}}{k^{1 - 1/n}} < \infty\,,
\end{equation}
then  almost all points on the curve\/ {\rm(\ref{bbkm:003})} are
not $\psi$-approximable. In the same paper Baker conjectured that
(\ref{bbkm:004}) could be replaced with the convergence of
(\ref{bbkm:002}), i.e.\ he conjectured that $\Veronese_n$ is of
Groshev type for convergence. This conjecture was proved by
V.I.~Bernik \cite{Bernik1} in 1989.

The divergence case was considered by V.V.~Beresnevich
\cite{Beresnevich9} in 1999 who proved that the Veronese curves\/
{\rm(\ref{bbkm:003})} are of Groshev type for divergence. The
proof is based on a new method involving regular systems,
introduced by Baker and Schmidt \cite{BakerSchmidt1} and used for
computing the Hausdorff dimension of sets of well approximable
points.

\subsection{Diophantine approximation on differentiable manifolds.}
\label{DA-on-manifolds}

In the sixties of the last century the investigations related to
the problem of Mahler eventually led to the development of a new
branch of metric number theory, usually referred to as
``Diophantine approximation of dependent quantities'' or
``Diophantine approximation on manifolds''. The first result
involving manifolds defined by functions satisfying some mild and
natural properties was obtained by Schmidt \cite{Schmidt3}, who
proved that any $C^{(3)}$ planar curve with curvature
non-vanishing almost everywhere is extremal. Schmidt's theorem was
subsequently improved by R.~Baker \cite{Baker5}, who has shown
that almost all points on Schmidt's curves are not
$\psi$-approximable whenever $(\ref{bbkm:004})_{n=2}$ is
satisfied. It has been recently shown that Schmidt's curves are of
Groshev type for convergence \cite{BernikDickinsonDodson1} and for
divergence \cite{BeresnevichBernikDickinsonDodson1}.

Until the mid-nineties most of the results in metric Diophantine
approximation dealt with manifolds of a special structure or of
high enough dimension. M.M.~Dodson, B.P.~Rynne and
J.A.G.~Vickers \cite{DodsonRynneVickers5,DodsonRynneVickers6,%
DodsonRynneVickers8} investigated a class of manifolds satisfying
a geometric condition which for surfaces in $\RR^3$ assumed two
convexity (e.g.\ a cylinder does not satisfy that condition).
Schmidt \cite{Schmidt3} has investigated certain straight lines in
$\RR^n$ for extremality, and recently such lines have been shown
to be of Groshev type \cite{BeresnevichBernikDickinsonDodson2}.

A new method, based on combinatorics of the space of lattices, was
developed in \cite{KleinbockMargulis1} by D.Y.~Kleinbock and
G.A.~Margulis\footnote{See also \cite{KleinbockMargulis2} and
\cite{Kleinbock1} for more on interactions between dynamics on the
space of lattices and metric Diophantine approximation.}, who
proved the extremality of the so-called non-degenerate manifolds
(also they proved these manifolds to be strongly extremal, see
Section~\ref{conc}).

\begin{definition}\label{dfn1}
Let $\vf:\U\longrightarrow\RR^n$ be a map defined on  an open set
$\U\subset\RR^d$. We say that $\vf$ is {\em $l$-non-degenerate
at}\/ $\vx_0\in\U$ if $\vf$ is $l$ times continuously
differentiable on some sufficiently small ball centered at $\vx_0$
and partial derivatives of $\vf$ at $\vx_0$ of orders up to $l$
span $\RR^n$. We say that $\vf$ is {\em non-degenerate at}\/
$\vx_0$ if it is $l$-non-degenerate at $\vx_0$ for some $l\in\NN$.
We say that $\vf$ is {\em non-degenerate}\/ if it is
non-degenerate almost everywhere on $\U$.
\end{definition}

The non-degeneracy of a manifold is naturally defined via the
non-degeneracy of its appropriate parameterization. Geometrically
the $l$-non-degeneracy of a manifold $\M\subset\RR^n$ at a point
$\vy_0\in\M$ means that for any hyperplane $\Pi$ in $\RR^n$,
$
   \limsup_{\vy\to\vy_0,\vy\in\M}\dist(\vy,\Pi)\cdot\|\vy-\vy_0\|^{-l}>0;
$
  that is,
the manifold can not be approximated by a hyperplane ``too well''
(see \cite{Beresnevich16,Beresnevich6}).

Recently Beresnevich \cite{Beresnevich16} (also a short version
published in \cite{Beresnevich12,Beresnevich13}), and
independently Bernik, Kleinbock and Margulis
\cite{BernikKleinbockMargulis2} using different techniques, have
proved that any non-degenerate manifold is of Groshev type for
convergence (also there is a multiplicative analogue and a more
general version of the result in \cite{BernikKleinbockMargulis2},
see Section~\ref{conc}).

Non-degenerate curves have been proved to be of Groshev type for
divergence \cite{Beresnevich14} (also
\cite{Beresnevich12,Beresnevich11} contain auxiliary parts of the
proof). Moreover, by Pyartli's method \cite{Pyartli1} one can
extend this result to analytic non-degenerate manifolds. The goal
of the present paper is to show that any non-degenerate manifold
is of Groshev type for divergence. The proof makes use of a new
technique, which involves a multidimensional analogue of regular
systems and extends the ideas of \cite{Beresnevich9}.

\subsection{The main result and the structure of the paper}
\label{result}

\begin{theorem}\label{maintheorem}
Let $\U$ be an open subset of $\RR^d$\/ and let\/
$\vf:\U\longrightarrow \RR^n$ be a non-degenerate map. Also let\/
$\psi:\RR_+\longrightarrow\RR_+$\/ be a non-increasing function
such that the sum {\rm (\ref{bbkm:002})} diverges. Then for almost
all $\vx\in\U$ the point $\vf(\vx)$ is $\psi$-approximable, i.e.\
for almost all $\vx\in\U$ there are infinitely many solutions
$\va\in\ZZ^n$ to the inequality
\begin{equation}\label{bbkm:005}
   |\langle\vf(\vx)\cdot\va\rangle|<\psi(\|\va\|_\infty^n).
\end{equation}
\end{theorem}

The proof of Theorem~\ref{maintheorem} is based on a method of
regular systems first suggested in \cite{Beresnevich9} for
dimension one. In particular, we generalize it for any dimension.
In Section~\ref{regular} we construct a regular system of resonant
sets corresponding to a given non-degenerate map. In
Section~\ref{approximation} we prove a general theorem on
approximation by resonant sets. And finally, Section~\ref{proof}
will complete the proof of Theorem~\ref{maintheorem}.

\section{Effective upper bounds}
\label{effective}

The result of this section will be applied to construct a regular
system of resonant sets. We show the following

\begin{theorem}\label{thm1}
Let $\vf:\U\longrightarrow\RR^n$ be non-degenerate at $\vx_0\in\U$. Then
there exists a sufficiently small ball $\B_0\subset\U$ centered at
$\vx_0$ and a constant $C_0>0$ such that for any ball
$\B\subset\B_0$ and any $\ve>0$ for all sufficiently big $Q$, one
has
\begin{equation}\label{bbkm:006}
   |\L_{\vf}(\B;\ve;Q)|\le C_0\,\ve\,|\B|,
\end{equation}
where
\begin{equation}\label{bbkm:007}
   \L_{\vf}(\B;\ve;Q)=\bigcup_{\va\in\ZZ^n\,:\ 0<\|\va\|_\infty\le Q}
   \set{\vx\in\B:|\langle\vf(\vx)\cdot\va\rangle|<\ve Q^{-n}}.
\end{equation}
\end{theorem}

The proof of Theorem~\ref{thm1} will rely on considering two
special cases: when the norm of the gradient $\va\nabla\vf(x)$ is
big, or, respectively, not very big. Theorem~\ref{thm2} below is
essentially due to Bernik and for $d=1$ has appeared earlier
\cite{Beresnevich14}. Its proof relies on the ideas of the method
of essential and inessential domains developed by Sprind\^zuk,
when he solved the problem of Mahler. Theorem~\ref{thm3} below is
due to Kleinbock and Margulis \cite{BernikKleinbockMargulis2} and
is proved by means of the method involving lattices, which was
first developed in \cite{KleinbockMargulis1}. The dichotomy of
big/small derivatives has been extensively used; in particular, it
is used in
\cite{Beresnevich12,Beresnevich16,BernikKleinbockMargulis2} to
prove the convergence case.

\begin{theorem}[Theorem~1.3 in \cite{BernikKleinbockMargulis2}]
\label{thm2} Let $\B_0 \subset\RR^d$ be a ball, and let $\vf\in
C^{(2)}(3\B_0)$. Fix $\delta>0$ and define
\begin{equation}\label{bbkm:008}
L_1 = \max_{\|\beta\|_1 = 2}\ \max_{\vx\in 2\B_0}\|\partial_\beta
\vf(\vx)\|_\infty.
\end{equation}
Then for every ball $\B\subset\B_0$ and any $\va\in\ZZ^n$ such
that
\begin{equation}\label{bbkm:009}
\|\va\|_\infty \ge \frac1{nL_1(\diam\B)^2}\,,
\end{equation}
the set
\begin{equation}\label{bbkm:010}
\L_\vf^{(1)}(\B;\delta;\va)=\Big\{\vx\in\B\,:\,\left\{
\begin{array}{c}
   |\langle \vf(\vx)\cdot\va \rangle|  < \delta, \\[1ex]
   \|\va\nabla\vf(\vx)\|_\infty \ge \sqrt{ndL_1\|\va\|_\infty}
\end{array} \ \ \Big\}
\right.
\end{equation}
has measure at most\/ $C_1\delta|\B|$, where $C_1>0$ is a constant
depending on $d$ only.
\end{theorem}

\begin{theorem}[Theorem~1.4 in \cite{BernikKleinbockMargulis2}]
\label{thm3} Let $\U\subset\RR^d$ be an open set, $\vx_0\in \U$,
and let $\vf:\U\longrightarrow\RR^n$ be a map $l$-non-degenerate at
$\vx_0$. Then there exists a ball $\B_0\subset\U$ centered at
$\vx_0$ such that $3\B_0\subset U$ with the following
property\/{\rm:} there exist a constant $C_2>1$ such that for any
ball $\B\subset\B_0$, any $\ve$ with $0<\ve<1$ and any $Q\ge 1$
the set
\begin{equation}\label{bbkm:011}
\L_{\vf}^{(2)}(\B;\ve;Q)=\bigcup_{\va\in\ZZ^n\,:\,0<\|\va\|_\infty\le
Q}\Big\{\vx\in \B\,:\, \left\{
\begin{array}{l}
   |\langle
\vf(\vx) \cdot\va \rangle|<\ve Q^{-n}, \\[1ex]
   \|\va\nabla
\vf(\vx) \|_\infty  < \sqrt{ndL_1Q}
\end{array}
\right.
  \ \Big\}
\end{equation}
satisfies
\begin{equation}\label{bbkm:012}
  |\L_{\vf}^{(2)}(\B;\ve;Q)|\le C_2(\ve
  Q^{-1/2})^{\textstyle\frac{1}{d(n+1)(2l-1)}}\cdot|\B|,
\end{equation}
where $L_1$ is defined in\/ {\rm(\ref{bbkm:008})}.
\end{theorem}

\begin{proof}[Proof of Theorem\/~{\rm\ref{thm1}}]

Fix a ball $\B_0$ as in the statement of Theorem~\ref{thm3} and
fix any ball $\B\subset\B_0$. It is easy to see that the set
$\L_{\vf}(\B;\ve;Q)$ is expressed as the following union of three
subsets
$$
\L_{\vf}(\B;\ve;Q)=
\left(\bigcup_{\va\in\ZZ^n\,:\,Q_1\le\|\va\|_\infty\le Q}
\L_{\vf}^{(1)}(\B;\ve Q^{-n};\va)\right)\bigcup
$$
\begin{equation}\label{bbkm:013}
\bigcup\L_{\vf}^{(2)}(\B;\ve;Q)\bigcup
\left(\bigcup_{\va\in\ZZ^n\,:\,\|\va\|_\infty\le Q_1}
\L_{\vf}^{(1)}(\B;\ve Q^{-n};\va)\right),
\end{equation}
where $Q_1=[1/(nL_1(\diam\B)^2)]+1$. The measure of the first
subset is estimated by Theorem~\ref{thm2}:
\begin{equation}\label{bbkm:014}
\left|\bigcup_{\va\in\ZZ^n\,:\,Q_1\le\|\va\|_\infty\le Q}
\L_{\vf}^{(1)}(\B;\ve Q^{-n};\va)\right|\le C_1\ve
Q^{-n}|\B|(2Q+1)^n.
\end{equation}

Next, for every $\va\in\ZZ^n$ such that $0<\|\va\|_\infty<Q_1$ we
obviously have
\begin{equation}\label{bbkm:015}
\L_{\vf}^{(1)}(\B;\ve Q^{-n};\va)\subset\L_{\vf}^{(1)}(\B;\ve
Q^{-n}Q_1;\va_1),
\end{equation}
where $\va_1=Q_1\va$. It is clear that $\|\va_1\|_\infty\ge Q_1$.
Therefore, we can apply Theorem~\ref{thm2} to the set in the right
hand side of (\ref{bbkm:015}). Thus,
\begin{equation}\label{bbkm:016}
   |\L_{\vf}^{(1)}(\B;\ve Q^{-n};\va)|\le |\L_{\vf}^{(1)}(\B;\ve
Q^{-n}Q_1;\va_1)|\le C_1\ve Q^{-n}Q_1|\B|
\end{equation}
Since the number of points $\va\in\ZZ^n$ with $0<\|\va\|_\infty\le
Q_1$ is less than $(2Q_1+1)^n$, we get
\begin{equation}\label{bbkm:017}
   \left|\bigcup_{\va\in\ZZ^n\,:\,\|\va\|_\infty\le Q_1}
\L_{\vf}^{(1)}(\B;\ve Q^{-n};\va)\right|\le (2Q_1+1)^nC_1\ve
Q^{-n}Q_1|\B|.
\end{equation}

On combining (\ref{bbkm:012}), (\ref{bbkm:014}), (\ref{bbkm:017})
and (\ref{bbkm:013}) and letting $C_0>2^nC_1$, we obtain
(\ref{bbkm:006}) for all sufficiently big $Q$. This completes the
proof of Theorem~\ref{thm1}.

\end{proof}

We will also use the following

\begin{lemma}[Lemma~6 in \cite{Beresnevich16}]
\label{qwer} Let $\alpha,\beta\in\RR_+$, $d\in\NN$, $\B$ be a ball
in $\RR^d$, $f:\B\longrightarrow\RR$ be a function such that $f\in
C^{(k)}$ and for some $j$ with $1\le j\le d$ one has
\begin{equation}\label{bbkm:018}
\inf_{\vv x\in \B}|\partial_j^{k}f(\vv x)|\ge \beta.
\end{equation}
Then
$$
\left|\{\vx\in \B:|f(\vv x)|\le\alpha\}\right|\le
3^{(k+1)/2}(k(k+1)/2+1)(\diam\B)^{d-1}\left(\frac{\alpha}{\beta}\right)^{1/k}.
$$
\end{lemma}

\section{Regular systems of resonant sets}
\label{regular}

\begin{definition}\label{defn1}
Let $\U$ be an open subset of $\RR^d$, $\R$ be a family of subsets
of $\RR^d$, $N:\R\longrightarrow\RR_+$ be a function and let $s$ be a
number satisfying $0\le s<d$. The triple $(\R,N,s)$ is called a
{\em regular system} in $\U$ if there exist constants
$K_1,K_2,K_3>0$ and a function $\lambda:\RR_+\longrightarrow\RR_+$ with
$\lim_{x\to+\infty}\lambda(x)=+\infty$ such that for any ball
$\B\subset\U$ and for any $T>T_0$, where $T_0=T_0(\R,N,s,\B)$ is a
sufficiently large number, there exist sets
\begin{equation}\label{bbkm:019}
   R_1,\dots,R_t\in\R \mbox{ \ with \ } \lambda(T)\le N(R_i)\le T
\mbox{ \ for \ }
   i=\overline{1,t}
\end{equation}
and disjoint balls
\begin{equation}\label{bbkm:020}
  \B_1,\dots,\B_t \mbox{ \ with \ } 2\B_i\subset\B \mbox{ \ for \ }
i=\overline{1,t}
\end{equation}
such that
\begin{equation}\label{bbkm:021}
  \diam(\B_i)=T^{-1} \mbox{ \ for \ } i=\overline{1,t},
\end{equation}
\begin{equation}\label{bbkm:022}
   t\ge K_1|\B|T^d
\end{equation}
and such that for any $\gamma\in\RR$ with $0<\gamma<T^{-1}$ one
has
\begin{equation}\label{bbkm:023}
K_2\gamma^{d-s}T^{-s}\le|\B(R_i,\gamma)\cap
\B_i|,\phantom{K_2\gamma^{d-s}T^{-s}\le}
\end{equation}
\begin{equation}\label{bbkm:024}
\phantom{\le K_3\gamma^{d-s}T^{-s},}|\B(R_i,\gamma)\cap 2\B_i|\le
K_3\gamma^{d-s}T^{-s}.
\end{equation}
The elements of $\R$ will be called {\em resonant sets}.
\end{definition}

This definition generalizes the concept of {\em regular system of
points}\/ of Baker and Schmidt. In fact, it is equivalent to the
Baker--Schmidt definition when $\U=\RR$, $\R$ consists of points
in the real line, and $s=0$ \cite{BakerSchmidt1}. In this
situation conditions (\ref{bbkm:023}) and (\ref{bbkm:024}) hold
automatically. Also this definition covers the multidimensional
concept of a regular system of points
\cite{Beresnevich11} when $s=0$.
Definition~\ref{defn1} is closely related to {\em ubiquitous
systems} \cite{DodsonRynneVickers3}.

The goal of this section is to establish the following

\begin{theorem}\label{thm4}
Let $\vf=(f_1,\dots,f_n):\U\longrightarrow\RR^n$ be a non-degenerate map,
where $\U$ is an open subset of $\RR^d$. Given an $\va\in\ZZ^n$,
$\va\not=0$ and an $a_0\in\ZZ$, let
$$
R_{\va,a_0}=\{\vx\in\U:\va\cdot\vf(\vx)+a_0=0\}.
$$
Define the following set
$$
\R_\vf=\{R_{\va,a_0}:\va\in\ZZ^n,\ \va\not=0,a_0\in\ZZ\}
$$
and the following function
$$
N(R_{\va,a_0})=(\|\va \|_\infty)^{n+1}.
$$
Then for almost every point $\vx_0\in\U$ there is a ball
$\B_0\subset\U$ centered at $\vx_0$ such that $(\R,N,d-1)$ is a
regular system in $\B_0$.
\end{theorem}

\begin{proof}
There is no loss of generality in assuming that $f_1(\vx)=x_1$. In
fact, using the non-degeneracy of $\vf$, it is possible to show
that $\vf'(\vx)\not=\vv0$ almost everywhere (see
\cite[Section~5]{Beresnevich16}). Thus we can take a sufficiently
small neighborhood of a point $\vx_0$ with $\vf'(\vx_0)\not=\vv0$
instead of the original domain $U$, and then make $f_1(\vx)$ equal
$x_1$ by an appropriate change of variables. Also, as $\vf$ is
non-degenerate, we can take $U$ to be a sufficiently small
neighborhood of a point $\vx_0$ such that $\vf$ is non-degenerate
at this point. Moreover,  we  can take $\B_0$ satisfying
Theorem~\ref{thm1}. Thus, in view of that theorem, for any ball
$\B\subset\B_0$ the set
$$
\G(\B;(4C_0)^{-1};Q)={\textstyle\frac34}\B
\smallsetminus\L_\vf({\textstyle\frac34}\B;(4C_0)^{-1};Q)
$$
will satisfy the estimate
\begin{equation}\label{bbkm:025}
|\G(\B;(4C_0)^{-1};Q)|\ge\frac12|\B|
\end{equation}
for all sufficiently large $Q$.

Note also that there is no loss of generality in assuming that
\begin{equation}\label{bbkm:026}
   \max_{1\le j\le n}\sup_{\vx\in\B_0}\|\nabla f_j(\vx)\|_\infty\le L_2,
\end{equation}
for some constant $L_2>0$.

The proof of Theorem~\ref{thm4} will be completed with the help of

\begin{proposition}\label{ppp1}
There is a sufficiently big number $Q_0$ such that for any $Q\ge
Q_0$ for any $\vx\in\G(\B;(4C_0)^{-1};Q)$ there is an integer
point $\va\in\ZZ^n$, $\va\neq\vv0$ and an integer $a_0$ with
\begin{equation}\label{bbkm:027}
Q^{n+1}=T/C_3\le N(R_{\va,a_0})\le T=C_3Q^{n+1},
\end{equation}
where $C_3=\Big(4C_0(nL_2)^{n-1}\Big)^{n+1}$, and a point $\vz\in
R_{\va,a_0}$ such that
\begin{equation}\label{bbkm:028}
\|\vx-\vz\|<C_4T^{-1},
\end{equation}
where $C_4=C_3n/(2C_0)$, and such that for any $\gamma$ with
$0<\gamma<T^{-1}$ we have
\begin{equation}\label{bbkm:029}
K_2\gamma T^{-(d-1)}\le |\B(R_{\va,a_0},\gamma)\cap
\B(\vz,T^{-1}/2)|,\phantom{K_2\gamma T^{-(d-1)}\le}
\end{equation}
\begin{equation}\label{bbkm:030}
\phantom{\le K_3\gamma T^{-(d-1)},}
|\B(R_{\va,a_0},\gamma)\cap\B(\vz,T^{-1})|\le K_3\gamma
T^{-(d-1)},
\end{equation}
where $K_2,K_3>0$ are some constants which depend on neither $\B$
nor $T$.
\end{proposition}

\begin{proof}[Proof of Proposition~{\rm\ref{ppp1}}]

Let $\vx\in\G(\B;(4C_0)^{-1};Q)$. By Minkowski's linear forms
theorem, there are integers $\va\in\ZZ^n$, $\va\neq\vv0$ and
$a_0\in\ZZ$ such that
\begin{equation}\label{bbkm:031}
\left\{
\begin{array}{l}
  |\vf(\vx)\cdot \va+a_0|\le (4C_0)^{-1}Q^{-n},\\[1ex]
  |a_1|\le 4C_0(nL_2)^{n-1}Q,\\[1ex]
  |a_i|\le Q/(nL_2) \ \ \ i=\overline{2,n}.
\end{array}
\right.
\end{equation}

Define the function $F(\vx)=\vf(\vx)\cdot\va+a_0$. It follows from
(\ref{bbkm:031}) that
\begin{equation}\label{bbkm:032}
\|\va\|_\infty\le 4C_0(nL_2)^{n-1}Q=T^{1/(n+1)}.
\end{equation}

Since $\vx\in\G(\B;(4C_0)^{-1};Q)$, $\|\va\|_\infty$ must be $>Q$,
which, combined with (\ref{bbkm:032}), gives (\ref{bbkm:027}).

As $|a_j|<Q$ for $j=\overline{2,n}$, we have $|a_1|>Q$. Now, using
(\ref{bbkm:026}) and the condition $f_1(\vx)=x_1$, we get
\begin{equation}\label{bbkm:033}
|\partial_1F(\vx)|=|a_1|\cdot|\partial_1f_1(\vx)|-
\sum_{i=2}^n|a_i|\cdot|\partial_1f_j(\vx)|>Q-
\sum_{i=2}^nQ/(nL_2)\cdot L_2=\frac{Q}{n}.
\end{equation}

Since $\partial_1\vf$ is uniformly continuous on $\B_0$, there is
a sufficiently small number $r_1>0$ such that for any
$\vx_1,\vx_2\in\U$ with $\|\vx_1-\vx_2\|<r_1$ we have
$$
\|\partial_1\vf(\vx_1)-\partial_1\vf(\vx_2)\|_\infty<\frac{1}{8n^2C_0(nL_2)^{n-1}}\,.
$$
It follows that
$$
|\partial_1 F(\vx_1)-\partial_1F(\vx_2)|\le
n\|\va\|_\infty\|\partial_1\vf(\vx_1)-\partial_1\vf(\vx_2)\|_\infty\le
\frac{1}{8nC_0(nL_2)^{n-1}}\,\|\va\|_\infty.
$$
Applying (\ref{bbkm:032}) now gives $|\partial_1
F(\vx_1)-\partial_1F(\vx_2)|\le Q/(2n)$ for all $\vx_1,\vx_2\in\U$
with $\|\vx_1-\vx_2\|<r_1$. This and (\ref{bbkm:033}) imply
$$
|\partial_1 F(\vy)|\ge |\partial_1
F(\vx)|-|\partial_1F(\vx)-\partial_1F(\vy)|>Q/(2n)
$$
for all $\vy\in\U$ with $\|\vx-\vy\|<r_1$.

As $\vx\in\frac34\B$, we have $\B(\vx,\diam\B/8)\subset\B$. Define
$r_0=\min(r_1,\diam\B/8)$. Thus,
\begin{equation}\label{bbkm:034}
|\partial_1 F(\vy)|>Q/(2n)\mbox{ \ for all $\vy\in\B(\vx,r_0)$}.
\end{equation}

\medskip

Let $|\theta|<r_0$. Then
$\vx_\theta=(x_1+\theta,x_2,\dots,x_d)\in\B(\vx,r_0)$, where
$\vx=(x_1,\dots,x_d)$. By the Mean Value Theorem, we have
$F(\vx_\theta)=F(\vx)+\partial_1F(\tvv x_\theta)\theta$, where
$\tvv x_\theta\in\B(\vx,r_0)$. This can equivalently be written as
\begin{equation}\label{bbkm:035}
  \frac{F(\vx_\theta)}{\partial_1F(\tvv
x_\theta)}=\frac{F(\vx)}{\partial_1F(\tvv x_\theta)}+\theta.
\end{equation}

Assume that $Q>(n/(2r_0C_0))^{1/(n+1)}$. This condition implies
that for any
$$
\theta\in\Big[-n/(2C_0)\cdot Q^{-n-1}\,,\,n/(2C_0)\cdot
Q^{-n-1}\Big]
$$
we have $|\theta|<r_0$, and therefore $\vx_{\theta},\tvv
x_{\theta}\in\B(\vx,r_0)$. Now using (\ref{bbkm:031}) and
(\ref{bbkm:034}) we get
$$
|F(\vx)/\partial_1F(\tvv x_\theta)|<n/(2C_0)\cdot Q^{-n-1}.
$$
It follows from this and (\ref{bbkm:035}) that
$F(\vx_\theta)/\partial_1F(\tvv x_\theta)$ is positive at
$\theta=n/(2C_0)\cdot Q^{-n-1}$ and negative at
$\theta=-n/(2C_0)\cdot Q^{-n-1}$. By continuity, there is a number
$\theta_0$ with
$$
|\theta_0|<n/(2C_0)\cdot Q^{-n-1}
$$
such that $F(\vx_{\theta_0})/\partial_1F(\tvv x_{\theta_0})=0$,
or, equivalently, $F(\vx_{\theta_0})=0$. Define $\vz$ to be
$\vx_{\theta_0}$. By construction, $\vz\in R_{\va,a_0}$, and
$\|\vx-\vz\|=|\theta_0|<n/(2C_0)\cdot Q^{-n-1}$. This proves
(\ref{bbkm:028}).

\medskip

{\it Now we are going to show}\/ (\ref{bbkm:030}). Assume that
$T>(C_4+1)/r_0$. This condition and (\ref{bbkm:028}) imply that
$$
\B(\vz,T^{-1})\subset\B(\vx,r_0).
$$

Let $0<\gamma<T^{-1}$. By definition, for any point
$\vy\in\B(R_{\va,a_0},\gamma)$ there is a point $\vy_0\in
R_{\va,a_0}$ such that $\|\vy-\vy_0\|<\gamma$.

Assume that $\vy\not=\vy_0$. Then, by the Mean Value Theorem, we
have
$$
F(\vy)=F(\vy_0)+\nabla F(\vy_1)\cdot(\vy-\vy_0)=\nabla
F(\vy_1)\cdot(\vy-\vy_0)=(\va\nabla \vf(\vy_1))\cdot(\vy-\vy_0),
$$
where $\vy_1$ is a point between $\vy_0$ and $\vy$. Using
(\ref{bbkm:032}), we find that
$$
|F(\vy)|\le d\|\va\nabla
f(\vy_1)\|_\infty\cdot\|\vy-\vy_0\|_\infty\le dn\|\va\|_\infty
L_2\gamma\le C_5Q\gamma,
$$
where $C_5=dn4C_0(nL_2)^{n-1}L_2$. It follows that
$$
\B(R_{\va,a_0},\gamma)\cap\B(\vz,T^{-1})\subset\{\vy\in
\B(\vz,T^{-1}):|F(\vy)|\le C_5Q\gamma\}.
$$
Now using Lemma~\ref{qwer}, this inclusion, (\ref{bbkm:034}), and
the fact that $\B(\vz,T^{-1})\subset\B(\vx,r_0)$, we obtain
$$
|\B(R_{\va,a_0},\gamma)\cap\B(\vz,T^{-1})|\le 12nC_5\gamma
T^{-(d-1)}.
$$
This implies inequality (\ref{bbkm:030}) with $K_3=12nC_5$.

\medskip

{\it It remains to show}\/ (\ref{bbkm:029}). If $d=1$, then
(\ref{bbkm:029}) holds with $K_2=1/2$. Thus we assume that $d>1$.

Define the constant
$$
C_6=\min\left\{1/8\,,\,\frac{1}{16(d-1)n^2L_2C_3^{1/(n+1)}}\right\}.
$$

Let $\vz'=(z_2,\dots,z_d)$, where $\vz=(z_1,\dots,z_d)$. Fix any
point $\vy'=(y_2,\dots,y_d)\in\RR^d$ such that
$\|\vy'-\vz'\|<C_6T^{-1}$. Given $y_1\in\RR$, we define the point
$\vy=(y_1,\vy')=(y_1,y_2,\dots,y_d)$. If $|y_1-z_1|\le T^{-1}/8$
then
\begin{equation}\label{bbkm:036}
\|\vy-\vz\|=\sqrt{|y_1-z_1|^2+\|\vy'-\vz'\|^2}\le
|y_1-z_1|+\|\vy'-\vz'\|<T^{-1}/8+C_6T^{-1}=T^{-1}/4.
\end{equation}
It follows that $\vy\in\B(\vz,T^{-1}/4)$ whenever $|y_1-z_1|\le
T^{-1}/8$. By the Mean Value Theorem,
$$
F(\vy)=F(\vz)+\nabla F(\tvv y)\cdot(\vy-\vz),
$$
where $\tvv y\in\B(\vz,T^{-1}/4)$. Since $F(\vz)=0$, we obtain
\begin{equation}\label{bbkm:037}
F(\vy)/\partial_1F(\tvv y)=(y_1-z_1)+\sum_{i=2}^d \partial_iF(\tvv
y)/\partial_1F(\tvv y)\cdot(y_i-z_i).
\end{equation}
Using (\ref{bbkm:032}), (\ref{bbkm:034}) and the inequality
$\|\vz'-\vy'\|< C_6T^{-1}$, we find that
$$
\left|\sum_{i=2}^d \partial_iF(\tvv y)/\partial_1F(\tvv
y)\cdot(y_i-z_i)\right|<T^{-1}/8.
$$
Therefore, the expression on the right of (\ref{bbkm:037}) is
positive when $y_1-z_1=T^{-1}/8$ and is negative when
$y_1-z_1=-T^{-1}/8$. Thus, the function
$f(y_1)=F(\vy)/\partial_1F(\tvv y)$ has different signs at $\pm
T^{-1}/8$. By the continuity, there is a point
$y_1\in(-T^{-1}/8,T^{-1}/8)$ such that $f(y_1)=0$, or,
equivalently, $F(y_1,\dots,y_d)=0$.

Thus, we have proved that for any $\vy'$ with
$\|\vy'-\vz'\|<C_6T^{-1}$ there is a point $y_1(\vy')\in\RR$ such
that $\vy=(y_1(\vy'),\vy')\in R_{\va,a_0}\cap\B(\vz,T^{-1}/4)$. It
is now easy to see that for any $\theta\in\RR$ with $|\theta|\le
T^{-1}/4$ we have $(y_1(\vy')+\theta,\vy')\in\B(\vz,T^{-1}/2)$.
Thus, for any positive $\gamma$ with $\gamma<T^{-1}$ the set
$$
A(\gamma)=\Big\{(y_1(\vy')+\theta,\vy'):\|\vy'-\vz'\|<C_6T^{-1},|\theta|\le
\gamma/4\Big\}
$$
satisfies
\begin{equation}\label{bbkm:038}
A(\gamma)\subset \B(R_{\va,a_0},\gamma)\cap\B(\vz,T^{-1}/2).
\end{equation}
By the theorem of Fubini, it is easy to calculate that
$$
|A(\gamma)|=|\B_{d-1}(\vz',C_6T^{-1})|_{d-1}\cdot \gamma/2=
|\B_{d-1}(\vv 0,C_6)|_{d-1}/2 \cdot\gamma\cdot T^{-(d-1)}.
$$
Applying (\ref{bbkm:038}) now gives inequality (\ref{bbkm:029})
with $K_2=|\B_{d-1}(\vv 0,C_6)|_{d-1}/2$.
\end{proof}

\medskip

{\it Now we proceed to prove Theorem~{\rm\ref{thm4}}.}

Assume that $Q>Q_0$. Choose a collection
$$
(\va_1,a_{0,1},\vz_1),\dots,(\va_t,a_{0,t},\vz_t)\in
(\ZZ^n\smallsetminus\{\vv0\})\times \ZZ\times \B \text{ with }
\vz_i\in R_{\va_i,a_{0,i}}
$$
such that
\begin{equation}\label{bbkm:039}
Q^{n+1}=T/C_3\le N(R_{\va_i,a_{0,i}})\le T=C_3Q^{n+1}\ (1\le i\le
t)
\end{equation}
and such that for any $\gamma$ with $0<\gamma<T^{-1}$ we have
\begin{equation}\label{bbkm:040}
K_2\gamma T^{-(d-1)}\le |\B(R_{\va_i,a_{0,i}},\gamma)\cap
\B(\vz_i,T^{-1}/2)|\ (1\le i\le t)\phantom{K_2\gamma
T^{-(d-1)}\le}
\end{equation}
\begin{equation}\label{bbkm:041}
\phantom{\le K_3\gamma T^{-(d-1)},}
|\B(R_{\va_i,a_{0,i}},\gamma)\cap\B(\vz_i,T^{-1})|\le K_3\gamma
T^{-(d-1)}\ (1\le i\le t)
\end{equation}
\begin{equation}\label{bbkm:042}
\B(\vz_i,T^{-1}/2)\cap\B(\vz_j,T^{-1}/2)=\varnothing\text{ for all
different $i,j$ ($1\le i,j\le t$)}
\end{equation}
and the number $t$ is maximal possible.

By Proposition~\ref{ppp1}, for any point
$\vx\in\G(\B;(4C_0)^{-1};Q)$ there is a triple
$$
(\va,a_{0},\vz)\in (\ZZ^n\smallsetminus\{\vv0\})\times \ZZ\times
\B \text{ with }\vz\in R_{\va,a_0}
$$
satisfying (\ref{bbkm:027}) --- (\ref{bbkm:030}). By the
maximality of $t$ there is an index $i\in\{1,\dots,t\}$ such that
$$
\B(\vz_i,T^{-1}/2)\cap\B(\vz,T^{-1}/2)\not=\varnothing.
$$
It follows that $\|\vz-\vz_i\|<T^{-1}$. This inequality and
(\ref{bbkm:028}) imply that $\|\vx-\vz_i\|<(C_4+1)T^{-1}$.
Therefore,
$$
\G(\B;(4C_0)^{-1};Q)\subset\bigcup_{i=1}^t\B(\vz_i,(C_4+1)T^{-1}).
$$
By this inclusion and (\ref{bbkm:039}), we obtain
$$
|\B|/2\le|\G(\B;(4C_0)^{-1};Q)|\le t\cdot|\B(\vv 0,C_4+1)|T^{-d}.
$$
Therefore, $t\ge K_1|\B|T^d$ with $K_1=(2|\B(\vv 0,C_4+1)|)^{-1}$.

Let $\lambda(x)=x/C_3$. Now, setting $R_i=R_{\va_i,a_{0,i}}$ and
$\B_i=\B(\vz_i,T^{-1}/2)$ gives the required collections of
resonant sets and balls in the definition of regular system. This
completes the proof of Theorem~\ref{thm4}.
\end{proof}

\section{Approximation by resonant sets}
\label{approximation}

In this section we prove the following general result, which is an
extension of Theorem~2 in \cite{Beresnevich9}.

\begin{theorem}\label{thm5}
Let $\U$ be an open set in $\RR^d$, and let $(\R,N,s)$ be a
regular system in $\U$. Let $\Psi\,{:}\,\RR_+\longrightarrow\RR_+$ be a
non-increasing function such that the sum
\begin{equation}\label{bbkm:043}
\sum_{h=1}^\infty h^{d-s-1}\Psi^{d-s}(h)
\end{equation}
diverges. Then for almost all points $\vv x\in \U$ the inequality
\begin{equation}\label{bbkm:044}
\dist(\vx,R)<\Psi(N(R))
\end{equation}
has infinitely many solutions $R\in\R$.
\end{theorem}

\subsection{Auxiliary lemmas}

\begin{lemma}\label{lem1}
Let $E\subset\RR^d$ be a measurable set, and let $\U\subset\RR^d$
be an open subset. Assume that there is a constant $\delta>0$ such
that for any finite ball $\B\subset\U$ we have $|E\cap\B|\ge
\delta|\B|$. Then $E$ has full measure in $\U$, i.e.\
$|\U\smallsetminus E|=0$.
\end{lemma}

\begin{proof}
Let $\tilde E=\U\smallsetminus E$. As $\U\smallsetminus\tilde
E=\U\cap E$, for any ball $\B\subset \U$ we have
$|\B\smallsetminus \tilde E|\ge \delta|\B|$. Next, for any $\ve>0$
there is a cover of $\tilde E$ consisting of balls $\B_i$ such
that
$$
\sum_{i=1}^\infty|\B_i|-\ve\le|\tilde E|\le
\sum_{i=1}^\infty|\B_i|.
$$
Notice that the sets $\B_i\smallsetminus\tilde E$ and
$\B_i\cap\tilde E$ are disjoint and satisfy $\B_i=
(\B_i\smallsetminus\tilde E)\cup(\B_i\cap\tilde E)$. Then we get
$$
|\tilde E|\ge \sum_{i=1}^\infty|\B_i|-\ve=
\sum_{i=1}^\infty|\B_i\smallsetminus\tilde E|+
\sum_{i=1}^\infty|\B_i\cap\tilde E|-\ve\ge
$$
$$
\ge\delta\sum_{i=1}^\infty|\B_i|+\left|\bigcup_{i=1}^\infty
\B_i\cap\tilde E\right|-\ve\ge \delta|\tilde E|+|\tilde E|-\ve.
$$
Therefore, $|\tilde E| \le\ve/\delta\to0$ as $\ve\to0$. Hence,
$\tilde E$ is null and $E$ has full measure in $\U$.
\end{proof}

\begin{lemma}[Lemma~5, Chapter~1 in \cite{Sprindzuk2}]\label{lll1}
Let $E_i\subset\RR^d$ be a sequence of measurable sets, and let
the set $E$ consist of points $\vx$ belonging to infinitely many
$E_i$. If there is a sufficiently large ball in $\RR^d$ which
contains all the sets $E_i$, and the sum $\sum_{i=1}^\infty|E_i|$
diverges, then
\begin{equation}\label{bbkm:045}
|E|\ge \limsup_{N\to\infty} \frac{\displaystyle\left(\sum_{i=1}^N
|E_i|\right)^2} {\displaystyle\sum_{i=1}^N\sum_{j=1}^N|E_i\cap
E_j|}\ .
\end{equation}
\end{lemma}

\begin{lemma}\label{lem10}
Let $\Psi$ satisfy the conditions of Theorem\/~{\rm\ref{thm5}},
and let $\tilde\Psi(h)=\min\{ch^{-1},\Psi(h)\}$, where $c>0$ is a
constant. Then $\tilde\Psi$ is non-increasing and the sum
\begin{equation}\label{bbkm:046}
\sum_{h=1}^\infty h^{d-s-1}\tilde\Psi^{d-s}(h)
\end{equation}
diverges.
\end{lemma}

\begin{proof}

The monotonicity of $\tilde\Psi$ is easily verified. Assume that
(\ref{bbkm:046}) converges. Then, by the monotonicity, we have
$$
l^{d-s}\tilde\Psi^{d-s}(l)\ll\sum_{l/2\le h\le
l}h^{d-s-1}\tilde\Psi^{d-s}(h)\to0\mbox{ as }l\to\infty.
$$
It follows that $l\tilde\Psi(l)=\min\{c,l\Psi(l)\}\to0$ as
$l\to\infty$. This is possible only if $l\Psi(l)\to0$ as
$l\to\infty$. It follows that $\tilde\Psi(l)=\Psi(l)$ for all
sufficiently large $l$. Therefore, the sum (\ref{bbkm:043})
converges, contrary to the conditions of Lemma~\ref{lem10}.

\end{proof}

\begin{lemma}\label{lem11}
Let $\Psi:\RR_+\longrightarrow\RR_+$ be non-increasing. Fix any $d>0$.
Then the sums
$$
\sum_{h=1}^\infty h^{d-1}\Psi^d(h)\mbox{ \ and \ }
\sum_{k=0}^\infty 2^{kd}\Psi^d(2^k)
$$
converge or diverge simultaneously.
\end{lemma}

\begin{proof}
Using the monotonicity of $\Psi$ we get the following inequalities
$$
2^{(k+1)d}\Psi^d(2^{k+1})\ll \sum_{2^k\le
h<2^{k+1}}h^{d-1}\Psi^d(h)\ll 2^{kd}\Psi^d(2^k).
$$
Summing these over all $k\in\NN$ gives the required property.
\end{proof}


\subsection{Proof of Theorem~\ref{thm5}}

By Lemma~\ref{lem10}, there is no loss of generality in assuming
that for all $h>0$
\begin{equation}\label{bbkm:047}
\Psi(h)\le h^{-1}/2.
\end{equation}

Fix any ball $\B\subset \U$ and set $T=2^k$. By
Definition~\ref{dfn1}, there are constants $K_1,K_2,K_3>0$, which
do not depend on $\B$, and there is a sufficiently big number
$k_0$ satisfying the following properties: for any natural number
$k\ge k_0$ there are resonant sets $R_k^{(i)}\in\R$ ($1\le i\le
t_k$) and balls $\B_k^{(i)}$ with $2 \B_k^{(i)}\subset \B$ ($1\le
i\le t_k$) such that
\begin{equation}\label{bbkm:048}
\lambda(2^k)\le N(R_k^{(i)})\le 2^k \ \ (1\le i\le t_k),
\end{equation}
\begin{equation}\label{bbkm:049}
\diam \B_k^{(i)}=2^k \ \ (1\le i\le t_k),
\end{equation}
\begin{equation}\label{bbkm:050}
\B_k^{(i)}\cap \B_k^{(j)}=\varnothing  \ \ (1\le i,j\le t_k,\
i\not=j),
\end{equation}
\begin{equation}\label{bbkm:051}
K_2\gamma^{d-s}2^{-sk}\le|\B(R_k^{(i)},\gamma)\cap \B_k^{(i)}|,
\end{equation}
and
\begin{equation}\label{bbkm:052}
|\B(R_k^{(i)},\gamma)\cap 2 \B_k^{(i)}|\le K_3\gamma^{d-s}2^{-sk}
\end{equation}
for any $\gamma$, $0<\gamma<2^{-k}$,
\begin{equation}\label{bbkm:053}
K_1 2^{dk}|\B|\le t_k\le 2^{dk}|\B|.
\end{equation}

For every natural number $k\ge k_0$ and $i\in\{1,\dots,t_k\}$ we
define the sets
$$
E_k^{(i)}= \B(R_k^{(i)},\Psi(2^k))\cap \B_k^{(i)}
$$
and
\begin{equation}\label{bbkm:054}
E_k=\bigcup_{i=1}^{t_k}E_k^{(i)}.
\end{equation}
It follows from (\ref{bbkm:051}) and (\ref{bbkm:052}) that
\begin{equation}\label{bbkm:055}
K_2\Psi^{d-s}(2^k)2^{-sk}\le|E_k^{(i)}|\le
K_3\Psi^{d-s}(2^k)2^{-sk}.
\end{equation}
It follows from (\ref{bbkm:050}) that
\begin{equation}\label{bbkm:056}
E_k^{(i)}\cap E_k^{(j)}=\varnothing\ \mbox{ if $i\not=j$, $1\le
i,j\le t_k$}.
\end{equation}
Therefore, $|E_k|=\sum_{i=1}^{t_k}|E_k^{(i)}|$. Using
(\ref{bbkm:053}) and (\ref{bbkm:055}), we find that
$$
K_1K_2\Psi^{d-s}(2^k)2^{(d-s)k}|\B|\le |E_k|\le
2K_1K_3\Psi^{d-s}(2^k)2^{(d-s)k}|\B|.
$$
Let $\phi_k=2^{(d-s)k}\Psi^{d-s}(2^k)$. Then we have
\begin{equation}\label{bbkm:057}
K_1K_2|\B|\phi_k\le|E_k|\le 2K_1K_3|\B|\phi_k.
\end{equation}

Using the divergence of (\ref{bbkm:043}) and applying
Lemma~\ref{lem11}, we obtain
\begin{equation}\label{bbkm:058}
\sum_{k=1}^\infty\phi_k=\infty.
\end{equation}

It follows from (\ref{bbkm:057}) and (\ref{bbkm:058}) that
$\sum_{k=k_0}^\infty|E_k|=\infty$. Since $\B$ is bounded  and all
the sets $E_k$ are contained in $\B$, Lemma~\ref{lll1} can be
applied to the sequence $E_k$. We are now going to obtain
estimates for the numerator and the denominator in
(\ref{bbkm:045}).

When $K>k_0$, inequalities (\ref{bbkm:057}) imply that
\begin{equation}\label{bbkm:059}
\sum_{k=k_0}^K|E_k|\ge K_1K_2|\B|\sum_{k=k_0}^K\phi_k.
\end{equation}

Now we proceed to estimate the measure of $E_k\cap E_l$. Let
$k_0\le k< l\le K$, where $K>k_0$. Using (\ref{bbkm:054}), we can
write
$$
E_l\cap E_k^{(i)}=\bigcup_{j=1}^{t_l}E_l^{(j)}\cap E_k^{(i)}.
$$

By (\ref{bbkm:055}), we find that $|E^{(j)}_l\cap E^{(i)}_k|\le
K_3\Psi^{d-s}(2^l)2^{-sl}$. Hence,
\begin{equation}\label{bbkm:060}
|E_l\cap E_k^{(i)}|\le K_3\Psi^{d-s}(2^l)2^{-sl}\cdot q(l,k,i),
\end{equation}
where $q(l,k,i)$ is the number of different indices $j$ such that
$E_l^{(j)}\bigcap E_k^{(i)}\not=\varnothing$.

Now we will estimate $q(l,k,i)$. Using (\ref{bbkm:049}) and
(\ref{bbkm:050}), we get
$$
\left|\bigcup_{j=1,\dots,t_l\,:\,E_l^{(j)}\cap
E_k^{(i)}\not=\varnothing} \B_l^{(j)}\right|= |\B(\vv
0,2^{-l}/2)|\cdot q(l,k,i)=\hspace*{11ex}
$$
\begin{equation}\label{bbkm:061}
\hspace*{14ex}=|\B(\vv 0,1/2)|\cdot2^{-dl}q(l,k,i).
\end{equation}

Consider any ball $\B_l^{(j)}$ such that $E_l^{(j)}\cap
E_k^{(i)}\not=\varnothing$. Fix a point $\vv x\in E_l^{(j)}\cap
E_k^{(i)}$. By the definition of $E_k^{(i)}$, there is a point
$\vv z\in R_k^{(i)}$ such that
\begin{equation}\label{bbkm:062}
   \|\vv x-\vv z\|<\Psi(2^k).
\end{equation}

Next, since $\vv x\in E_l^{(j)}\subset\B^{(j)}_l$, for any point
$\vv y\in \B^{(j)}_l$ we have
\begin{equation}\label{bbkm:063}
\|\vv y-\vv x\|<\diam \B^{(j)}_l=2^{-l}.
\end{equation}
Then, using (\ref{bbkm:062}) and (\ref{bbkm:063}), we obtain
$$
\|\vv y-\vv z\|<\|\vv y-\vv x\|+\|\vv x-\vv z\|<2^{-l}+\Psi(2^k).
$$
Therefore,
$$
  \dist(\vv y,R^{(i)}_k)<2^{-l}+\Psi(2^k),
$$
whence
\begin{equation}\label{bbkm:064}
\B^{(j)}_l\subset \B(R^{(i)}_k,2^{-l}+\Psi(2^k)).
\end{equation}

Let $\vv x_0$ denote the center of $\B^{(i)}_k$. For $\vv x\in
E_k^{(i)}\subset\B^{(i)}_k$, we have
\begin{equation}\label{bbkm:065}
   \|\vv x-\vv x_0\|<\frac12\diam\B^{(i)}_k.
\end{equation}
Using the inequality $l>k$ and (\ref{bbkm:063}), we obtain
$$
\|\vv x-\vv y\|<\frac12\diam\B^{(i)}_k.
$$
On combining the last inequality with (\ref{bbkm:065}), we get
$$
\|\vv y-\vv x_0\|<\|\vv y-\vv x\|+\|\vv x-\vv
x_0\|<\diam\B^{(i)}_k.
$$
Thus, $\B^{(j)}_l\subset2\B^{(i)}_k$. Using this inclusion and
(\ref{bbkm:064}) gives
\begin{equation}\label{bbkm:066}
\bigcup_{j=1,\dots,t_l\,:\,E_l^{(j)}\cap
E_k^{(i)}\not=\varnothing}
  \B_l^{(j)}\subset  \B(R^{(i)}_k,2^{-l}+\Psi(2^k))\cap2 \B^{(i)}_k.
\end{equation}
Now, applying (\ref{bbkm:052}), (\ref{bbkm:066}), and the
monotonicity of the measure, we derive
$$
\left|\bigcup_{j=1,\dots,t_l\,:\,E_l^{(j)}\cap
E_k^{(i)}\not=\varnothing} \B_l^{(j)}\right|\le
K_3(2^{-l}+\Psi(2^k))^{d-s}2^{-sk}\le \hspace*{10ex}
$$
$$
\hspace*{20ex}\le K_32^{d-s}(2^{-l(d-s)}+\Psi^{d-s}(2^k))2^{-sk}.
$$
On combining this inequality and (\ref{bbkm:061}), we obtain
\begin{equation}\label{bbkm:067}
q(l,k,i)\ll 2^{s(l-k)}+2^{dl}2^{-sk}\Psi^{d-s}(2^k).
\end{equation}
It follows from (\ref{bbkm:060}) and (\ref{bbkm:067}) that
\begin{equation}\label{bbkm:068}
|E_l\cap E_k^{(i)}|\ll \Psi^{d-s}(2^l)2^{-sk}+
2^{(d-s)l}2^{-sk}\Psi^{d-s}(2^k)\Psi^{d-s}(2^l).
\end{equation}

Since the number of different sets $E^{(i)}_k$ does not exceed
$t_k$, we have
$$
|E_l\cap E_k|\le t_k\cdot\max_{1\le i\le t_k}|E_l\cap E^{(i)}_k|.
$$
Using this inequality, (\ref{bbkm:053}), and (\ref{bbkm:068}), we
get
$$
|E_l\cap E_k|\ll |\B| \Psi^{d-s}(2^l)2^{(d-s)k}\times
$$
\begin{equation}\label{bbkm:069}
\times\Big(1+2^{(d-s)l}\Psi^{d-s}(2^k)
\Big)=|\B|(2^{-(d-s)(l-k)}\phi_l+\phi_k\phi_l).
\end{equation}

For arbitrary $l,k$ with $k_0\le l,k\le K$, we have
\begin{equation}\label{bbkm:070}
|E_l\cap E_k|\ll |\B|(2^{-(d-s)|l-k|}\phi_l+\phi_k\phi_l).
\end{equation}

By (\ref{bbkm:058}), there is a sufficiently big number $K'$ such
that for all $K>K'$
\begin{equation}\label{bbkm:071}
\sum_{k=k_0}^K\phi_k>1.
\end{equation}
Let $K>K'$. Now using (\ref{bbkm:057}), (\ref{bbkm:069}), and
(\ref{bbkm:071}), we calculate
$$
\sum_{l=k_0}^K\sum_{k=k_0}^K|E_l\cap E_k|\ll
|\B|\sum_{l=k_0}^K\sum_{k=k_0}^K\phi_k\phi_l+
|\B|\sum_{l=k_0}^K\sum_{k=k_0}^K2^{-(d-s)|l-k|}\phi_l\le
$$
$$
\le |\B|\sum_{l=k_0}^K\sum_{k=k_0}^K\phi_k\phi_l+
|\B|\sum_{l=k_0}^K\phi_l\sum_{k\in\ZZ}2^{-(d-s)|l-k|}=
$$
$$
=|\B|\left(\sum_{l=k_0}^K \phi_l\right)^2+
(1+2^{s-d})/(1-2^{s-d})|\B|\sum_{l=k_0}^K\phi_l\ll
|\B|\left(\sum_{k=k_0}^K\phi_k\right)^2
$$
where the implicit constant in this estimate does not depend on
either $\B$ or $K$. Using (\ref{bbkm:059}) now gives
\begin{equation}\label{bbkm:072}
\frac{\left(\sum_{k=k_0}^K\ |E_k|\right)^2}
{\sum_{l=k_0}^K\sum_{k=k_0}^K\ |E_l\cap E_k|}\gg |\B|
\end{equation}
when $K>K'$. By Lemma~\ref{lll1}, the set $E$ consisting of points
$\vv x$ which belong to infinitely many sets $E_k$ has measure
$\ge(K_1K_2)^2/C_{10}\cdot|\B|$.

Using the monotonicity of $\Psi$ and inequalities
(\ref{bbkm:048}), it is easy to see that for any point $\vv x \in
E$ inequality (\ref{bbkm:044}) has infinitely many solutions. Let
$\R(\Psi)$ denote the set of points $\vv x\in \U$ such that
inequality (\ref{bbkm:044}) has infinitely many solutions. Then
$$
E\subset \R(\Psi)\cap\B.
$$
It follows that
$$
|\R(\Psi)\cap\B|\ge |E|\gg |\B|.
$$
By Lemma~\ref{lem1}, the set $\R(\Psi)$ has full measure in $\U$.
The proof of Theorem~\ref{thm5} is completed.

\section{Proof of the main theorem}
\label{proof}

It is obvious that we can restrict ourselves to a sufficiently
small ball $\B_0$ centered at a point belonging to a set with full
measure in $U$. By Theorem~\ref{thm4} we can take $\B_0$ to be
such that $(\R,N,s)$ is a regular system in $\B_0$, where $s=d-1$,
$N$ and $\R$ are defined in the statement of Theorem~\ref{thm4}.
Define the sequence $\Psi$ by setting
$$
dnL_2 h\Psi(h^{n+1})=\psi(h^n).
$$
Thus $\Psi(k)=k^{-1/(n+1)}\psi(k^{n/(n+1)})/dnL_2$. Since $\psi$
is non-increasing, $\Psi$ is non-increasing as well. Next, we
calculate
$$
\sum_{h=1}^\infty h^{d-s-1}\Psi^{d-s}(h)= \sum_{h=1}^\infty
\Psi(h)=
$$
$$
=\frac{1}{dnL_2}\sum_{k=1}^\infty\sum_{(k-1)^{(n+1)/n}< h\le
k^{(n+1)/n}} h^{-1/(n+1)}\psi(h^{n/(n+1)})\gg
$$
$$
\gg\sum_{k=1}^\infty\sum_{(k-1)^{(n+1)/n}<h\le k^{(n+1)/n}}
k^{-1/n}\psi(k)\ge \sum_{k=1}^\infty \psi(k)=\infty.
$$

By Theorem~\ref{thm5}, for almost all $\vx\in\U$ there are
infinitely many $(\va,a_0)\in\ZZ^n\times\ZZ$ satisfying
\begin{equation}\label{bbkm:073}
\dist(\vx,R_{\va,a_0})<\Psi(\|\va\|_\infty^{n+1}).
\end{equation}

It follows from (\ref{bbkm:073}) that there is a point $\vz\in
R_{\va,a_0}$ such that
\begin{equation}\label{bbkm:074}
\|\vx-\vz\|<\Psi(\|\va\|_\infty^{n+1}).
\end{equation}

By the definition of $R_{\va,a_0}$, we have
$F(\vz)=\va\cdot\vf(\vz)+a_0=0$. Using the Mean Value Theorem, we
obtain
\begin{equation}\label{bbkm:075}
F(\vx)=F(\vz)+\nabla F(\tvv x)\cdot(\vx-\vz)= \nabla F(\tvv
x)\cdot(\vx-\vz)=(\va\nabla \vf(\tvv x))\cdot(\vx-\vz).
\end{equation}
where $\tvv x$ is a point between $\vx$ and $\vz$. Using
(\ref{bbkm:026}), we find that
\begin{equation}\label{bbkm:076}
|\,\langle\va\cdot\vf(\vx)\rangle\,|=|F(\vx)|\le d\|\va\nabla
f(\tvv x)\|_\infty\cdot\|\vx-\vz\|_\infty<dn\|\va\|_\infty L_2
\Psi(\|\va\|_\infty^{n+1})=\psi(\|\va\|_\infty^n).
\end{equation}

As we have shown above, for almost all $\vx\in\U$ there are
infinitely many $(\va,a_0)\in\ZZ^n\times\ZZ$ satisfying
(\ref{bbkm:073}). Therefore, for almost all $\vx\in\U$ there are
infinitely many $\va$ satisfying (\ref{bbkm:076}). This completes
the proof of Theorem~\ref{maintheorem}.

\section{Concluding remarks}\label{conc}

In this section we give a brief account of other results in metric
Diophantine approximation and state the most important problems in
this field. Also we discuss possible developments of the theory of
regular systems and difficulties that prevent us from proving
multiplicative divergence Khintchine type results.

\subsection{Simultaneous approximation}

The point $\vy\in\RR^n$ is called {\it simultaneously
$\psi$-approxi\-mable}\/ if
\begin{equation}\label{bbkm:077}
\|\,\langle q\vy\rangle\,\|_\infty^n<\psi(q)
\end{equation}
has infinitely many solutions $q\in\ZZ$. By the Khintchine
transference principle, a point $\vy\in\RR^n$ is very well
approximable if and only if it is simultaneously
$\psi_\ve$-approximable for some positive $\ve$, where
$\psi_\ve(h)= h^{-(1+\ve)}$. Unfortunately there is no such
connection between simultaneous and dual approximation for general
approximation functions $\psi$ that would make it possible to
derive a Khintchine type theorem for the simultaneous case from
the dual and visa verse. However, it has been known since the 1926
paper of Khintchine that almost all (almost no) points of $\RR^n$
are simultaneously $\psi$-approximable if the sum (\ref{bbkm:002})
diverges (converges).

Let $\M$ be a submanifold of $\RR^n$. One says that $\M$ is of
{\em Khintchine type for divergence $($for convergence$)$}\/ if
almost all (almost no) points of $\M$ are simultaneously
$\psi$-approximable whenever the sum (\ref{bbkm:002}) diverges
(converges).

We mostly deal with monotonic approximation errors. However, it is
worth saying that for $n>1$ an analogue of Khintchine's theorem
for non-monotonic error function has been obtained by
A.~Pollington and R.~Vaughan \cite{PollingtonVaughan1}, who proved
a multidimensional analogue of the Duffin--Schaeffer conjecture.

Only special manifolds have been proved to be of Khintchine type.
Bernik \cite{Bernik3} has shown that the parabola
$\{(x,x^2):x\in\RR\}$ is of Khintchine type for convergence. He
has also proved with a method of trigonometric sums that any
manifold given as a topological product of at least 4 planar
curves with curvatures non-vanishing almost everywhere is of
Khintchine type for both convergence and divergence
\cite{Bernik4}. A class of manifolds in $\RR^n$ with a special
geometrical property, which substantially restricts the dimension
of the manifolds, has been proved to be of Khintchine type for
both convergence and divergence
\cite{DodsonRynneVickers6,DodsonRynneVickers8}.

In the Khintchine type theory for simultaneous Diophantine
approximation the following is regarded as the main problem.

\Problem\label{pr3} Prove that a non-degenerate manifold $\M$ in
$\RR^n$ is of Khintchine type for convergence and for divergence.

It is of interest to consider some special cases of
Problem~\ref{pr3} such as the circle, the sphere and others. There
remain two classical special cases of Problem~\ref{pr3}: to prove
that for $n\ge3$ the curve $\Veronese_n$ is of Khintchine type for
convergence and to prove that for $n\ge2$ the curve $\Veronese_n$
is of Khintchine type for divergence.

One difficulty in the simultaneous Diophantine approximation is
that there is no longer the dichotomy of big/small derivative (the
derivative is always big) but the investigated sets are quite
rare. Thus one needs a considerably new technique to break through
the problem.

A much deeper problem is to prove asymptotic formulae for the
number of solutions of Diophantine inequalities under
consideration. This remains unsettled for both linear and
simultaneous approximation.

\subsection{Multiplicative results}
The point $\vv y\in\RR^n$ is said to be {\em
$\psi$-multiplicatively approximable}\/ if the inequality
\begin{equation}\label{bbkm:078}
|\,\langle\va\cdot\vv x\rangle\,|<\psi(\Pi_+(\va))
\end{equation}
has infinitely many solutions $\va\in\ZZ^n$, where
$
\Pi_{+}(\va)=\prod_{i = 1}^n \max(|a_i|,1).
$
One can define {\em very well multiplicatively approximable}\/
points to be $\psi_\ve$-multiplicatively approximable for some
positive $\ve$, with $\psi_\ve(h)=h^{-1-\ve}$.

By the Borel--Cantelli lemma, almost all points of $\RR^n$ are not
$\psi$-multiplicatively approximable whenever the sum
\begin{equation}\label{bbkm:079}
\sum_{h=1}^{\infty}(\log h)^{n-1}{\psi(h)}
\end{equation}
converges. Since $\Pi_{+}(\va)$ is not greater than
$\|\va\|_\infty^n$, any $\psi$-approximable point is automatically
$\psi$-multiplicatively approximable. Therefore, a very well
approximable point is also very well multiplicatively
approximable.

A manifold $\M$ is said to be of {\em multiplicative Groshev type
for divergence $($convergence$)$}\/ if almost all (almost no)
points of $\M$ are multiplicatively $\psi$-approximable whenever
the sum (\ref{bbkm:079}) diverges (converges). A manifold $\M$ is
said to be {\em strongly extremal}\/ if almost all points of $\M$
are not very well multiplicatively approximable.

The problem of proving strong extremality in connection with
multiplicative approximation was first raised by Baker in
\cite[Ch.~9, p.~96]{Bak90}. The question, as initially proposed,
related to the Veronese curve and it was later generalized to any
non-degenerate manifold by Sprindzuk. Baker was motivated in part
by the non-metrical instances of specific points known to have the
property of strong extremality, {\em i.e.}\/ the algebraic numbers
and powers of $e$ \cite[Ch.~7~and~Ch.~10]{Bak90}.

Kleinbock and Margulis \cite{KleinbockMargulis1} proved that any
non-degenerate manifold is strongly extremal, and later jointly
with Bernik \cite{BernikKleinbockMargulis1} they have shown a
stronger result that these manifolds are of multiplicative Groshev
type for convergence. They even proved a more general result, to
be stated in Section~\ref{conc_gen}. No manifold (except $\RR^n$
itself) has ever been shown to be of multiplicative Groshev type
for divergence.

The difficulty of proving multiplicative Groshev type theorems for
divergence with the method of this paper is that Minkowski's
theorem on convex bodies cannot be efficiently extended to
non-convex bodies, e.g.\ star bodies, which appear in the context
of multiplicative approximation. One might try to relax the
definition of regular system used in this paper by taking a
multi-valued function $N$ to control any possible difference in
the magnitude of integer coefficients. But in this way one would
loose a sufficient estimate for denominators in (\ref{bbkm:045}).
Thus more investigation is required to prove a multiplicative
Groshev type theorem for divergence.

\Problem\label{pr4} Prove that any non-degenerate manifold is of
multiplicative Groshev type for divergence.

One can also consider a multiplicative version of simultaneous
Diophantine approximation when one replaces the right hand side of
(\ref{bbkm:077}) with $\prod_{i = 1}^n |\langle qy_i\rangle|$.
Khintchine type theorems for this type of approximation have never
been proved for convergence or for divergence.

\subsection{A general approximation function}\label{conc_gen}

Let $\Psi:\ZZ\longrightarrow\RR_+$, $n,m\in\NN$. The point $\vv
y\in\RR^{nm}$ is said to be {\em $(\Psi,n,m)$-approximable}\/ if
the inequality
\begin{equation}\label{bbkm:080}
\|\,\langle\va\vv y\rangle\,\|_\infty^m<\Psi(\va)
\end{equation}
has infinitely many solutions $\va\in\ZZ^n$. The point $\vv y$ is
considered to be a matrix with $n$ rows and $m$ columns.

Due to Schmidt \cite{Schmidt5,Schmidt2} one knows the following
most general result on Diophantine approximation of independent
quantities.

Let $m,n\in\NN$, $n\ge2$, $\Psi:\ZZ^n\longrightarrow\RR_+$. Almost all
(almost no) points $\vv y\in\RR^{nm}$ are
$(\Psi,n,m)$-approximable whenever the sum
\begin{equation}\label{bbkm:081}
\sum\limits_{\va\in \ZZ^n} {\Psi(\va)}
\end{equation}
diverges (converges).

For the case of $m=1$ and under some monotonicity restrictions on
$\Psi$, Bernik, Kleinbock and Margulis extended the convergence
part of this result to non-degenerate manifolds. More precisely,
assuming that for every $i = \overline{1,n}$
\begin{equation}\label{bbkm:082}
\Psi(q_1,\dots,q_i,\dots,q_n) \ge \Psi(q_1,\dots,q'_i,\dots,q_n)
\quad\text{whenever} \quad |q_i| \le |q'_i|\text{ and }q_iq'_i >
0,
\end{equation}
they proved that almost no point $\vv y\in\M$ is
$(\Psi,n,1)$-approximable whenever the sum
(\ref{bbkm:081})${}_{m=1}$ converges, where $\M$ is a given
non-degenerate manifold.

\Problem\label{pr5} Assuming (\ref{bbkm:082}), prove that almost
all points $\vv y\in\M$ are $(\Psi,n,1)$-approximable whenever the
sum (\ref{bbkm:081}) diverges, where $\M$ is a given
non-degenerate manifold.

It is also of interest to investigate Diophantine approximation
(of any type) with non-monotonic error function (right hand side
of inequalities).

Another interesting problem is to find reasonable conditions of
the entries of the matrix $\vy$ in (\ref{bbkm:080}) when they are
dependent, so that one would have an extremality type or
Khintchine--Groshev (or Schmidt) type theorem.

\subsection{Hausdorff dimension}

The first results on the Hausdorff dimension of sets arising in
Diophantine approximation are due to V.~Jarnik and
A.S.~Besicovitch.
  They found
the exact value of the Hausdorff dimension of the set of
$w$-approximable points (i.e.\ $\psi_{w/n-1}$-approximable points
with $\psi_\ve(h)=h^{-1-\ve}$) in the real line.

The first general method for obtaining lower bounds for the
Hausdorff dimension was suggested by Baker and Schmidt. They
introduced the concept of regular systems, which made it possible
to efficiently describe the distribution of objects that were used
for approximation. Baker and Schmidt have proved with their method
that the set of $w$-approximable points on $\Veronese_n$ has
dimension at least $\frac{n+1}{w+1}$, and conjectured that this
number is the right upper bound as well. The Baker--Schmidt
conjecture was proved by Bernik \cite{Bernik6} in 1983. Extending
the ideas of Baker and Schmidt, Dodson and H.~Dickinson
\cite{DickinsonDodson1} have shown that for any extremal manifold
$\M$ in $\RR^n$ the set of $w$-approximable points on the manifold
has Hausdorff dimension at least $\frac{n+1}{w+1}+\dim\M-1$. Thus
we've got a very natural

\Problem\label{pr6} Let $w>n$ and $\M$ be a non-degenerate
manifold in $\RR^n$. Prove that the Hausdorff dimension of
$w$-approximable points on $\M$ is exactly
$\frac{n+1}{w+1}+\dim\M-1$.

Also, Dodson \cite{Dodson1,Dodson2} has investigated the Hausdorff
dimension of the set of $(\Psi,n,m)$-approximable points when
$\Psi(\va)=\psi(\|\va\|_\infty)$, $\psi:\RR_+\longrightarrow\RR_+$ and
Dickinson and S.~Velani \cite{DickinsonVelani1} answered a very
general question of the Hausdorff measure (with respect to a
general dimension function) of this set and proved a
Khintchine--Groshev type theorem.

The problem of calculating the Hausdorff dimension in the case of
simultaneous Diophantine approximation seems to be even more
complicated (see \cite[pp.~92--98]{BernikDodson1}). The Hausdorff
dimension of simultaneously $v$-approximable points (i.e.\
simultaneously $\psi_{(nv-1)}$-approximable with
$\psi_\ve(h)=h^{-1-\ve}$) with $v$ big enough seems to depend on
arithmetic and other properties of the manifold one would like to
approximate (see \cite[pp.~90--98]{BernikDodson1}). However there
might be a general formula for $v$ close to the extremal exponent
$1/n$.

\subsection{Beyond the non-degeneracy condition}
Looking for new classes of extremal, Khintchine or Groshev type
manifolds is a challenging task. The simplest ones for which the
non-degeneracy condition fails are proper affine subspaces of
$\RR^n$.  They have been studied in several papers in the past,
and some conditions (written in terms of Diophantine properties of
coefficients of parametrizing equations) have been found
sufficient for their extremality \cite{Schmidt3, Sprindzuk2} and,
in the case of straight lines passing through the origin,  for
being of Groshev type for both convergence and divergence
\cite{BeresnevichBernikDickinsonDodson2}.

Recently, in a preprint \cite{Kleinbock2} by Kleinbock, using the
dynamical approach of \cite{KleinbockMargulis1}, necessary and
sufficient conditions for extremality and strong extremality of
any  affine subspace of $\RR^n$ have been written down. Also it
has been shown there that a smooth submanifold $\M$ of an  affine
subspace $\L$ of $\RR^n$ is extremal (resp.~strongly  extremal)
whenever $\L$ is such, provided $\M$ is {\it non-degenerate in\/}
$\L$. The latter notion is a straightforward generalization of
Definition~\ref{dfn1}, so that a submanifold $\M$ of $\L$ is
non-degenerate in $\L$ if it can not be ``too well'' approximated
by hyperplanes contained in  $\L$.

This naturally leads to the following

\Problem\label{pr7} Find criteria for an  affine subspace $\L$ of
$\RR^n$ being of Groshev type for convergence or divergence; or,
given a specific function $\psi$ such that the sum
(\ref{bbkm:002}) diverges (converges), find necessary and
sufficient conditions for almost all (almost no) points of $\L$
being $\psi$-approximable. Also, prove that the aforementioned
properties of $\L$ are inherited by its submanifolds which are
non-degenerate in $\L$.

It is also worthwhile to mention that one can investigate
Diophantine properties of almost all (almost no) points with
respect to measures other than Lebesgue measures on smooth
manifolds. The latter can be supported on fractal subsets of
$\RR$ (see \cite{Weiss}) or $\RR^n$
(\cite{KleinbockLindenstraussWeiss}, the work currently in
progress).

\medskip

{\bf Acknowledgments.} Bernik, Kleinbock and Margulis are grateful
to SFB-343 and Humboldt Foundation for the support of their 1999
stay  at the University of Bielefeld, where a preliminary form of
this paper was discussed. Beresnevich is grateful to EPSRC for the
support of his stay at the University of York in 2000, where a
part of the paper was written, and to Prof.~Maurice~Dodson for his
hospitality and making arrangements for this stay.

\providecommand{\bysame}{\leavevmode\hbox
to3em{\hrulefill}\thinspace}

\end{document}